\documentstyle[12pt]{article}

\input amssym.def
\input amssym.tex

\def\R{\Bbb R}

\def\Z{\Bbb Z}

\def\C{\Bbb C}
\def\CP{\Bbb C\Bbb P}

\newtheorem{statement}{Statement}

\newenvironment{remark}
{\smallskip\noindent{\bf Remark\/}.}{\smallskip\par}

\newenvironment{proof}
{\noindent{\bf Proof\/}.}{{ $\Box$}\smallskip\par}

\title{On $A_k$--singularity on a plane curve of fixed degree.}
\author{S.M.Gusein--Zade
\thanks{Partially supported by grants 
RFBR--98--01--00612 and INTAS--96--0713. The work was done while the
author enjoyed the hospitality of the University of Nice.} 
\and N.N.Nekhoroshev
}
\date{}

\begin{document}

\def\eps{\varepsilon}

\maketitle

There is a general problem to describe singularities which can be met on
algebraic hypersurfaces, in particular on plane curves, of fixed degree
(see, e.g., \cite{GLS}). Here we shall consider $A_k$--singularities which
can be met on a plane curve of degree $d$. Let $k(d)$ be the maximal possible
integer $k$ such that there exists a plane curve of degree $d$ with an
$A_k$--singularity. Statement \ref{st1} gives an upper bound for $k(d)$.
According to it $\overline{\lim}_{d\rightarrow\infty}k(d)/d^2\le 3/4$.
We construct a
plane curve of degree $28s+9$ ($s\in\Z_{\ge 0}$) which has an
$A_k$--singularity with $k=420s^2+269s+42$. Therefore one has
$\underline{\lim}_{d\rightarrow\infty}k(d)/d^2\ge 15/28$ (pay attention that
$15/28>1/2$). The example is constructed basically in the same way as a
curve of degree $22$ with an $A_{257}$ singularity in \cite{GN} (the
aim of that example was somewhat different).

\begin{statement}{\label{st1}}
$k(d)\le(d-1)^2-\left[{d\over 2}\right]\cdot
\left(\left[{d\over 2}\right]-1\right).$
\end{statement}

\begin{remark}
We believe that this statement is known, however we have not found
a reference for it.
\end{remark}

\begin{proof}
Without any loss of generality one can suppose that the curve is reduced.
Therefore one can choose the infinite line in the complex projective plane
$\CP^2$ so that it intersects
the curve at $d$ different points. Let the affine part of the curve
$d$ be given by the equation $\{f(x, y)=0\}$ (where $f$ is a polynomial
of degree $d$). The surface $\{f(x,y)+z^2=\eps\}\subset\C^3$ is nonsingular
and the negative inertia index of the intersection form on its second
homology group $H_2(V; \R)$ just coincides with the right hand side of the
inequality in the statement (this follows, e.g., from the result of
J.Steenbrink on the intersection form for quasihomogeneous
singularities (\cite{S}), applied to the (isolated) singularity
$f_d(x,y)+z^2$, where $f_d(x,y)$ is the homogeneous part of degree $d$
of the polynomial $f$. Now the statement follows from the facts that
the intersection form of the $A_k$--singularity of $3$ variables is
negative defined and, if the considered curve has an $A_k$--singularity, 
then the vanishing homology group of this singularity is embedded into
the homology group $H_2(V; \R)$ (see, e.g., \cite{AGV}).
\end{proof}

\begin{statement}{\label{st2}}
There exists a plane curve of degree $28s+9$ which has an $A_k$--singularity
with $k=420s^2+269s+42$ ($s\in\Z_{\ge0}$).
\end{statement}

\begin{proof}
Let $\ell=3s+1$, $m=7s+2$, and let an (affine plane) curve $C$ be
given by the equation
\begin{equation}\label{eq1}
F(x, y)=y^2-2yA(x,y)+x^{8\ell}+4x^{7\ell}y^{m}=0,
\end{equation}
where $A(x,y)=\left[x^{4\ell}+2x^{3\ell}y^{m}-2x^{2\ell}y^{2m}
+4x^{\ell}y^{3m}-10y^{4m}\right]$. 
The degree $d$ of the curve $C$ is equal to $4m+1=7\ell+m=28s+9$.
Let $z=y-A(x,y)$ ($x$ and $z$ are local coordinates on the plane
near the origin). Then
\begin{equation}\label{eq2}
F(x,y)=z^2+56x^{3\ell}y^{5m}-56x^{2\ell}y^{6m}
+80x^{\ell}y^{7m}-100y^{8m}.
\end{equation}
Since $z(x,y)=y-x^{4\ell}+$ terms of higher degree, one has 
$y(x,z)=z+x^{4\ell}+$ terms of higher degree. Substituting this 
expression into (\ref{eq2}) one gets
$$F(x,z)=z^2+ 56 x^{3\ell+20\ell m}+\sum a_{ij}x^iz^j=
z^2+ 56 x^{k+1}+\sum a_{ij}x^iz^j,
$$
where $k=420 s^2+269s+42$ and the sum contains only members with
powers $(i, j)$ which lie over the segment $(0,2)(k+1,0)$, i.e.,
those for which $i/(k+1)+j/2>1$. This proves the statement.
\end{proof}

\begin{remark}
There exists a similar problem to determine the maximal number
$k'=k'(d)$ such that the singularity $A_{k'}$ is an adjacent 
to an isolated homogeneous singularity of degree $d$ (or (what
is equivalent) to a singularity from the same stratum $\mu=const$).
One has $k'(d)\ge k(d)$. All the reasonings above are valid for this
problem as well and thus the inequalities of the statements \ref{st1}
and \ref{st2} hold also for the number $k'(d)$.
\end{remark}

\noindent Moscow State University,\newline
Dept. of Mathematics and Mechanics,\newline
Moscow, 119899, Russia.\newline
{E-mail:} sabir\symbol{'100}mccme.ru\newline
\hspace*{1pt}{\ \ \ \ \ \ \ \ \ \ }nekhoros\symbol{'100}nw.math.msu.su

\end{document}